\documentclass{article}

\usepackage{amsfonts}

\begin{document}

\newcommand\A{{\mathcal A}}
\newcommand\eps{{\epsilon}}
\newcommand\fol{{\text{Fol}}}
\newcommand\hcs{{h_{\mathcal C}}}
\newcommand\hc{{h^*_{\mathcal C}}}
\newcommand\hcSs{{h_{\mathcal S\mathcal C}}}
\newcommand\hcS{{h^*_{\mathcal S\mathcal C}}}
\newcommand\id{{\text{Id}}}
\newcommand\LL{{\mathcal L}}
\newcommand\MCs{{\mathcal C}}
\newcommand\MC{{\mathcal C}^*}
\newcommand\MD{{\mathcal D}}
\newcommand\MDs{{\mathcal D}^*}
\newcommand\MG{{\mathcal G}}
\newcommand\MS{{\Sigma(\mathcal G)}}
\newcommand\muh{{\hat{\mu}}}
\newcommand\new[1]{{\bf #1}}
\newcommand\N{{\mathbf N}}
\newcommand\QED{$\square$}
\newcommand\R{{\mathbf R}}
\newcommand\Sh{{\hat{\Sigma}}}
\newcommand\Shs{{\hat{\Sigma}}^*}
\newcommand\subsetneq{\subset_\ne}
\newcommand\supsetneq{\supset_\ne}
\newcommand\T{{\mathbf T}}
\newcommand\text[1]{{\hbox{\rm #1}}}
\renewcommand\top{{\text{top}}}
\newcommand\tr{{\text{Tr}}}
\newcommand\X{{\mathcal X}}
\newcommand\Z{{\mathbf Z}}

\newtheorem{coro}{Corollary}
\newtheorem{defi}{Definition}
\newtheorem{lem}{Lemma}
\newtheorem{prop}{Proposition}
\newtheorem{theo}{Theorem}

\newenvironment{demo}
      {\medbreak\noindent{\sc Proof:}}
      {\hfill\QED\medbreak}

\newenvironment{demof}[1]
      {\medbreak\noindent{\sc Proof of {#1}:}}
      {\hfill\QED\medbreak}

\newenvironment{remk}
      {\medbreak\noindent{\sl Remark.}}
      {\medbreak}

\newenvironment{exam}
      {\medbreak\noindent{\sl Example.}}
      {\medbreak}

\newenvironment{claim}
      {\medbreak\noindent{\bf Claim.} \sl}
      {\medbreak}

\begin{title}
{Subshifts of Quasi-Finite Type}
\end{title}

\begin{author}
{J\'er\^ome Buzzi\\
{\tt buzzi@math.polytechnique.fr}
$ $\\
Centre de Math\'ematiques U.M.R. 7640\\
C.N.R.S. \& Ecole polytechnique}
\end{author}

\begin{date}
{March 2003}
\end{date}

\maketitle

\begin{abstract}
{We introduce subshifts of quasi-finite type as a generalization of the
well-known subshifts of finite type. This generalization is much less rigid
and therefore contains the symbolic dynamics of many non-uniform
systems, e.g., piecewise monotonic maps of the interval with positive entropy. 
Yet many properties remain: existence of finitely many ergodic
invariant probabilities of maximum entropy; lots of periodic points;
meromorphic extension of the Artin-Mazur zeta function.}
\end{abstract}

\section{Introduction}

Jacques Hadamard \cite{Hadamard} founded symbolic dynamics in 1898
when he realized that the dynamics of the geodesic flow on
surfaces of negative curvature can be represented by very simple
subsets of $\A^\Z$ ($\A$ being some finite subset). Namely, these
subsets are defined by excluding a finite number of words. Such
subsets are now called \new{subshifts of finite type} (or S.F.T.).
They have been thoroughly studied (see, e.g., \cite{LM}) and the
result of Hadamard has been generalized to all {\sl uniformly
hyperbolic systems} (see, e.g., \cite{Shub}). However, S.F.T. are
much too rigid to provide a description of more general dynamics
(for instance, there are only countably many topological conjugacy
classes of S.F.T.).

Therefore a key problem is to enlarge S.F.T. to accomodate wide
classes of {\bf non-uniform} dynamics and yet keep
most of the basic features of S.F.T.

\medbreak

In this paper, we provide a solution by introducing a new class of
subshifts, which we call \new{subshifts of quasi-finite type}.
They include the symbolic dynamics of a large class of
non-uniform dynamical systems: piecewise monotonic
maps \cite{dMvS} with positive entropy and more generally
entropy-expanding maps \cite{EE} satisfying a technical
assumption.

Whereas subshifts of finite type are described by finitely many
constraints, we allow a slowly growing number of constraints of a
given length, "slow growth" meaning with a rate strictly less than
the topological entropy.

We prove that these subshifts of quasi-finite type remarkably have
the same basic properties as S.F.T. at least with respect to
"complexity":
 \begin{itemize}
  \item they have finitely many ergodic
invariant probability measures maximizing entropy;
  \item they have lots of periodic points;
 \item their Artin-Mazur zeta functions have meromorphic extensions.
\end{itemize}
\medbreak

This paper can be considered as yet another illustration of the
following principle \cite{EE}:  {\sl  complexity
bounds imply semi-uniform hyperbolicity.}

\subsection{Definitions}

We consider a \new{subshift}, i.e., a closed $\sigma$-invariant
subset $\Sigma\subset\A^\Z$ ($\A$ is some finite set, $\sigma$
denotes the left-shift on $\A^\Z$). The one-sided version of
$(\Sigma,\sigma)$ is $(\Sigma_+,\sigma_+)$ with
$\Sigma_+:=\{A_0A_1A_2\dots:A\in\Sigma\}\subset\A^\Z$ and
$\sigma_+$ the left-shift on these one-sided sequences.

It is customary to consider \new{follower sets} \cite{LM}: if $A_{-n}\dots
A_0$ is some finite word on the alphabet $\A$, then
 $$
    \fol(A_{-n}\dots A_0) := \{B_0B_1B_2\dots : B\in\Sigma \text{
    and } B_{-n}\dots B_0 = A_{-n}\dots A_0\} \subset [A_0]_+ \subset \Sigma_+.
 $$
By convention, the follower of the empty word is the whole of
$\Sigma_+$. We can also write:
 $$
   \fol(A_{-n}\dots A_0) = \sigma_+^n([A_{-n}\dots A_0]_+) \subset \Sigma_+
 $$
where $[A_{-n}\dots A_0]_+$ denotes the cylinder in $\Sigma_+$.

\begin{defi}
A \new{left constraint} is a finite word $A_{-n}\dots A_0$,
$n\geq0$, such that\footnote{Observe that $\fol(A_{-n}\dots
A_0)=\emptyset$ if and only if the given word does not appear in
$\Sigma$.}:
 $$
      \emptyset \ne \fol(A_{-n}\dots A_0)\subsetneq \fol(A_{-n+1}\dots
              A_0).
 $$
The set of left constraints with length $n$ is denoted by
$\MCs(\Sigma,n)$ (or just $\MCs(n)$).

The \new{left constraint entropy} is the quantity:
 $$
      \hcs(\Sigma):=\limsup_{n\to\infty} \frac1n\log^+ \#\MCs(\Sigma,n).
 $$
\end{defi}

\begin{defi}
The \new{symmetric constraint entropy} is:
 $$
    \hcSs(\Sigma) = \min(\hcs(\Sigma),\hcs(\bar\Sigma)).
 $$
where $\bar\Sigma:=\{(A_n)_{n\in\Z}:(A_{-n})_{n\in\Z}\in\Sigma\}$.
\end{defi}

\medbreak

 Recall that the \new{topological entropy} is:
 $$
  h_\top(\Sigma)=\limsup_{n\to\infty} \frac1n\log \#\LL(\Sigma,n)
 $$
where $\LL(\Sigma,n):=\{A_0\dots A_{n-1}\in\A^n: [A_0\dots
A_{n-1}]_\Sigma\ne\emptyset\}$ where $[A_0\dots A_{n-1}]_\Sigma$
(or simply $[A_0\dots A_{n-1}]$) denotes a cylinder in $\Sigma$.

\medbreak

 We are at least in position to define the main object of this
paper:

\begin{defi}
$\Sigma$ is a \new{subshift of quasi-finite type} (or Q.F.T.) iff:
 $$
      \hcSs(\Sigma)<h_\top(\Sigma).
 $$
\end{defi}

\subsection{Examples and relatives}

\begin{lem} \label{lem-sofic}
All S.F.T. and sofic shifts with non-zero topological entropy are
Q.F.T. More precisely, $\hcs(\Sigma)=0$ if $\Sigma$ is a S.F.T.
(but $\hcs(\Sigma)>0$ is possible for sofic $\Sigma$).
\end{lem}

On the other hand, many symbolic dynamics which are not S.F.T. or
sofic are Q.F.T.:

\medbreak

A \new{piecewise monotonic map} \cite{dMvS} is a map $f:I\to I$ on some compact interval $I$
such that there is a finite partition of $I$ into subintervals on each of
which the restriction of $f$ is continuous and strictly monotonic. The natural
partition $P$ is the collection of maximum open intervals on which $f$ is
continuous and strictly monotonic. The symbolic dynamics is:
 $$
   \Sigma(f):= \{ A\in P^\Z:\forall n\in\Z\forall k\geq 0\;
                       <A_n\dots A_{n+k}>\ne\emptyset \}
 $$
where the notation $<A_0\dots A_k>:=A_0\cap f^{-1}A_1\cap\dots\cap
f^{-k}A_k\subset I$ stands for the \new{geometric cylinders}.

\medbreak

 A \new{multi-dimensional $\beta$-transformation}
\cite{Affine} is a map $T:[0,1)^d\to[0,1)^d$, with $d\geq1$ and
$T(x)=B.x\text{ mod }\Z^d$ where $B$ is an
expanding\footnote{i.e., $\exists\lambda>1\; \forall x,y\in\R^d\;
\|B.x-B.y\|\geq\lambda\|x-y\|$.} affine map of $\R^d$. The natural
partition $P$ is the finite collection of maximum open subsets of
$(0,1)^d$ on which $T(x)-B.x$ is constant. The symbolic dynamics
is defined as above. According to \cite{Affine}, this is a special
case of {\it connected piecewise entropy-expanding map}:

\medbreak

A \new{piecewise entropy-expanding map} is $(X,P,f)$ with (see
\cite{MEM}):
 \begin{itemize}
  \item $X$ is a compact subset of some Euclidean space;
  \item $P$ is a finite collection of pairwise disjoint open
  subset of $X$;
  \item $f:\bigcup_{A\in P} P\to X$ is such that each restriction
  $f:A\to f(A)$ can be extended to a homeomorphism between
  neighborhoods of $\bar A$ and $\overline{f(A)}$;
  \item the fundamental inequality:
     $$
        h_B(X,P,f):=\limsup_{n\to\infty} \frac1n\log \#\left\{A\in P^n: <A>\cap
        \bigcup_{B\in P} \partial f(B) \ne\emptyset\right\} < h_\top(\Sigma(f))
     $$
 \end{itemize}
Its symbolic dynamics $\Sigma(f)$ is again defined in the same way.

Notice that entropy-expanding does not imply expanding.

$(X,P,f)$ is said to be \new{connected} if every $P$-cylinder is
 connected.

\begin{lem} \label{lem-examples-qft}
The symbolic dynamics of the following dynamical systems are
Q.F.T.:
\begin{enumerate}
    \item piecewise monotonic maps with positive topological
entropy.
    \item connected piecewise entropy-expanding maps (hence, in particular,
multi-dimensional $\beta$-transformations).
\end{enumerate}
\end{lem}

This implies immediately:

\begin{coro} \label{coro-entropy}
The entropy of a Q.F.T. can take any value in $(0,\infty)$. In
particular, there are uncountably many conjugacy classes of
Q.F.T., in contrast to the case of S.F.T.
\end{coro}

There is an important weakening of Q.F.T.:

\begin{defi}
A left constraint $A_{-n}\dots A_0$, $n\geq0$, is
\new{extendable} if there exist a sequence $B\in\A^{-\Z}$ with $B_{-n}
\dots B_0= A_{-n}\dots A_0$ and infinitely many integers $m\geq0$
such that: $B_{-m}\dots B_0$ is again a left
constraint.\footnote{Remark that this is strictly stronger than
requiring that $A_{-n}\dots A_0$ is a suffix of infinitely many
left constraints.}

The set of extendable left constraints with length $n$ is denoted
by $\MC(\Sigma,n)$ (or just $\MC(n)$).

The \new{extendable left constraint entropy} is
 $$
    \hc(\Sigma) := \limsup_{n\to\infty}
    \frac1n\log\#\MC(\Sigma,n).
 $$
Subshifts $\Sigma$ with
$\hcS(\Sigma):=\min(\hc(\Sigma),\hc(\bar\Sigma))<h_\top(\Sigma)$
are called \new{weak-Q.F.T.}.
\end{defi}

Weak-Q.F.T. are definitely not as nice as Q.F.T.:

\begin{lem}\label{lem-hc-bad}
There exists a weak-Q.F.T. with countably infinitely many maximum
measures.
\end{lem}

The following qualititative generalization of sofic shifts is a
special case of weak-Q.F.T.:

\begin{defi} \label{def-e-M}
$\Sigma$ is \new{eventually Markovian on the left} iff for each
$A\in\Sigma$ there exists an integer $N$ such that:
 $$
    \forall n\geq N\quad \fol(A_{-n}\dots A_0) = \fol(A_{-N}\dots
    A_0).
 $$
\end{defi}

\begin{lem} \label{lem-e-M}
If $\Sigma$ is eventually Markovian on the left, then
$\hc(\Sigma)=0$.

In particular, all sofic shifts have $\hc(\Sigma)=0$ (compare with
$\hcs(\Sigma)$) and are weak-Q.F.T. (if they have positive
topological entropy).
\end{lem}

We exhibit some facts that show that the refinements of our
definitions (extendability condition, symmetry) do enlarge the
class of subshifts under consideration.

\begin{lem} \label{lem-R}
There are subshifts $\Sigma$ such that
$\hcs(\bar\Sigma)\ne\hcs(\Sigma)$ and even such that:
 $$
    \hcs(\bar\Sigma) < \hcs(\Sigma)=h_\top(\Sigma).
 $$

The same is true for $\hc$.
\end{lem}

\begin{lem}\label{lem-hc-hcs}
Obviously, we have:
 $$
    0\leq \hc(\Sigma) \leq\hcs(\Sigma) \leq h_\top(\Sigma)
 $$
and there exist subshifts for which all these inequalities are
strict. In particular, there are weak-Q.F.T. which are not Q.F.T.
\end{lem}

We now compare Q.F.T. to previously studied notions:

\medbreak

\new{Coded systems with synchronizing words} of Blanchard and Hansel
\cite{Blanchard} have a significant intersection with weak Q.F.T.:

\begin{lem} \label{lem-good-coded}
Any coded system such that the set of sequences not containing a
synchronizing word has topological entropy $<h_\top(\Sigma)$ is a
weak-Q.F.T.
\end{lem}

All coded systems are topologically transitive\footnote{i.e.,
there exists $A\in\Sigma$ such that $\{\sigma^nA:n\geq0\}$ is
dense in $\Sigma$.}, in contrast to Q.F.T. Hence we have trivial
examples of Q.F.T. which are not coded. Anne Bertrand
\cite{Bertrand} characterized (one-dimensional)
$\beta$-transformation with symbolic dynamics which are coded
systems with synchronizing words. In particular, not all of them
have this property. Hence:

\begin{lem} \label{lem-not-coded}
There are topologically transitive Q.F.T. which are not coded
systems with synchronizing words.
\end{lem}

We note the condition introduced by B.M. Gurevich \cite{Gurevich2}
which involves explicitly the speed with which finite order Markov
topological chains approximate the subshift. It seems unrelated to
Q.F.T.

\medbreak

We have the following relationships of $\hc(\Sigma)$ and
$\hcs(\Sigma)$ with known characteristics:

\begin{lem} \label{lem-K-cap}
The \new{boundary capacity} defined by Keller \cite{Keller}:
 $$
     \text{cap}(\Sigma) :=
        \limsup_{n\to\infty} \frac1n\log
        \sup_{r\geq1,\; w\in\LL(\Sigma,r)}
        \#\{ a\in\LL(\Sigma,n):
            [a]_+\cap\sigma_+^r[w]\ne\emptyset\qquad\qquad\qquad\qquad
 $$
 $$
   \qquad\qquad\qquad\qquad\qquad\qquad\qquad\qquad\qquad\text{ and }
            [a]_+\cap(\Sigma_+\setminus\sigma_+^r[w])\ne\emptyset \}
 $$
satisfies $\text{cap}(\Sigma)\geq\hcs(\Sigma)$ and the inequality
may be strict.
\end{lem}

\begin{lem} \label{lem-forbid}
The \new{entropy of minimal forbidden words} $h_{\mathcal
M}(\Sigma)$ considered by B\'eal and others \cite{Beal}:
 $$
   h_{\mathcal M}(\Sigma) := \limsup_{n\to\infty} \frac1n
      \#\{A_1\dots A_n: [A_1\dots A_n]=\emptyset \text{ and }
           [A_1\dots A_{n-1}]\ne\emptyset \qquad\qquad\qquad\qquad
 $$
 $$
 \qquad\qquad\qquad\qquad\qquad\qquad\qquad\qquad\qquad\qquad\qquad\qquad \text{ and }
           [A_2\dots A_n]\ne\emptyset \}
 $$
satisfies $h_{\mathcal M}(\Sigma)\leq\hcs(\Sigma)$, the inequality
being strict in some cases.

On the other hand, both $h_{\mathcal M}(\Sigma)> \hc(\Sigma)$ and
$h_{\mathcal M}(\Sigma)< \hc(\Sigma)$ occur. Moreover, there
exists a subshift $\Sigma$ with $h_{\mathcal M}(\Sigma)=0$ with
$\hc(\Sigma)=h_\top(\Sigma)>0$.\footnote{I would be interested by
an example with $h_{\mathcal M}(\Sigma)<h_\top(\Sigma)$ and
infinitely many maximal measures.}
\end{lem}

Finally we turn to the natural-looking \new{follower entropy}:
 $$
     h_\fol(\Sigma) := \limsup_{n\to\infty} \frac1n\log\#\{
        \fol(w) : w\in \LL(\Sigma,n) \}.
 $$
It is natural to ask whether this gives rise to a reasonable
variant of Q.F.T. This is not the case:

\begin{lem} \label{lem-hF-bad}
There exist subshifts with $h_\fol(\Sigma)=0<h_\top(\Sigma)$ which
have uncountably many \new{maximum measures}.\footnote{i.e.,
ergodic invariant probability measures with maximum entropy.}
\end{lem}

\subsection{Basic properties}

\begin{prop}\label{prop-top-invariance}
$\hcs(\Sigma)$ and $\hc(\Sigma)$ (and therefore their symmetric
versions) are invariants of topological conjugacy. On the other
hand, $\hcs(\Sigma)$ does not necessarily decrease under factor
maps also the Q.F.T. and weak-Q.F.T. properties are not preserved
under extensions or factors.
\end{prop}

\noindent {\bf Question:} {\sl Can $\hc(\Sigma)$ also increase
under factor maps?}

\begin{lem} \label{lem-constructions}
For unions:
 $$
  \hcs(\Sigma_1\cup\Sigma_2) = \max(\hcs(\Sigma_1),\hcs(\Sigma_2))
 $$
and for products:
 $$
  h_\top(X\times Y)-\hcs(X\times Y) = \min(h_\top(X)-\hcs(X),h_\top(Y)-\hcs(Y))
$$
In particular, if $X$ and $Y$ are Q.F.T., then $X\times Y$ is also
a Q.F.T.

These properties are also true for $\hc$.
\end{lem}

In particular, a product of an arbitrary subshift with a subshift
of zero-entropy is never Q.F.T or even weak-Q.F.T.

\medbreak

We shall see the following topological properties:

\begin{lem} \label{lem-top}
A Q.F.T. is not necessarily topologically transitive. A
weak-Q.F.T. always contains periodic points hence it is never
topologically minimal.
\end{lem}

\section{Main results}

\begin{theo} \label{theo-main}
Let $\Sigma$ be a Q.F.T. Then,
 \begin{enumerate}
  \item $\Sigma$ admits a finite number of maximum measures, each one of which is
  Bernoulli (up to a period\footnote{$\mu$ is Bernoulli up to a period $p\geq1$ iff
  there exists a subset $X$ such that $\sigma^p(X)=X$, $\Sigma$ is the disjoint union
  $\bigcup_{k=0}^{p-1}\sigma^k X$ (up to a negligible subset) and $(\sigma^p,\mu|X)$
  is Bernoulli.});
  \item the Artin-Mazur zeta function:
     $$
       \zeta(z) := \exp - \sum_{n\geq1} \frac{z^n}{n}
            \#\{x\in\Sigma: \sigma^nx=x\}
     $$
 extends from a holomorphic function on $|z|<e^{-h_\top(\Sigma)}$ to a
meromorphic function on $|z|< e^{-\hcSs(\Sigma)}$. Moreover, $|z|=
e^{-\hcSs(\Sigma)}$ is the natural boundary\footnote{i.e., $\zeta$
cannot be extended meromorphically to a connected set
$U\supsetneq\{|z|< e^{-\hcSs(\Sigma)}\}$.} of meromorphic
extension for some Q.F.T.\footnote{This is known for the
symbolic dynamics of $\beta$-transformations for Lesbegue-a.e.
$\beta>1$.}
  \item the number of periodic points satisfies:
     $$
        0 < \limsup_{n\to\infty} \frac{ \#\{x\in\Sigma:\sigma^nx=x\} }
                { e^{n h_\top(\Sigma)}} < \infty.
     $$
  \end{enumerate}

If $\Sigma$ is only weak-Q.F.T. then the following properties remain true:
 \begin{enumerate}
 \item there are at most a {\em countable} number of
maximum measures (but there can be infinitely many of them --- see
Lemma \ref{lem-hc-bad});
 \item the zeta function may fail to have a meromorphic extension (see also Lemma
\ref{lem-hc-bad}).
 \item $\limsup_{n\to\infty}(1/n)\log\#\{x\in\Sigma:\sigma^nx=x\}=
 h_\top(\Sigma)$.
 \end{enumerate}
\end{theo}

We recall that a \new{Markov shift} is defined as follows. Given a
{\sl countable} oriented graph $\mathcal G$, the associated Markov
shift $\Sigma(\mathcal G)$ is the set of all paths on $\mathcal
G$:
 $$
    \Sigma(\mathcal G) := \{ g\in\mathcal G^\Z : \forall n\in\Z\;
       g_n \longrightarrow g_{n+1} \text{ in } \mathcal G \}
 $$
together with the left-shift $\sigma$. Observe that if we require
the graph to be finite, then this reduces to S.F.T.

 \medbreak

The theorem above will follow from the following structure
theorem:

\begin{theo} \label{theo-struc}
Let $\Sigma\subset\A^\Z$ be a weak-Q.F.T. Then there is a
countable oriented graph $\MG$  and a map $\pi$ from the set of
vertices of $\MG$ to $\A$ such that the induced map
$\pi:\MS\to\Sigma$:
 \begin{itemize}
    \item is well-defined and satisfies $\pi\circ\sigma=\sigma\circ\pi$;
    \item is one-to-one between $\MS$ and $\Sigma\setminus\mathcal X$ with:
     \begin{enumerate}
     \item $
          \limsup_{n\to\infty} \frac1n \log^+ \#\mathcal X\cap
                \{x:\sigma^n x=x\}  \leq \hcS(\Sigma)$;
     \item $h(\sigma,\mu)\leq \hcS(\Sigma)$ for all invariant
     probability measures $\mu$ with $\mu(\mathcal X)=1$.
     \end{enumerate}
 \end{itemize}

Given any $h>\hcSs(\Sigma)$, there are only finitely many
irreducible parts of $\MG$ with entropy larger than $h$.
\end{theo}

In the language of \cite{EE},

\begin{coro}
A weak-Q.F.T. is entropy-conjugate to a Markov shift.
\end{coro}

\begin{remk}
This theory can easily accomodate {\sl weights}, i.e., one can
introduce a potential function $\psi:\Sigma\to\R$ and define
constraint pressure $P_{\mathcal C}(\Sigma)$ and consider
subshifts with $P_{\mathcal C}(\Sigma)<P_\top(\Sigma)$, the usual
topological pressure w.r.t. the fixed potential. Then all the
above results hold (or rather their weighted counterparts), with
the possible exception of the extendability of the zeta-function
(the corresponding result for Markov shifts has not been proved as
far as I know).
\end{remk}

\subsection*{Some questions}

\begin{itemize}
  \item Could we prove the above theorem by direct methods, i.e.,
  without using Markov diagrams? (This is possible in a geometric
  setting by using induction on a Markov rectangle.)
  \item It would be more elegant to have a single construction instead of
breaking the left/right symmetry.
  \item The above theorem ensures the existence of a
``good'' presentation for any Q.F.T. On the other hand,  what is
the set of presentations of QFT?
  \item Does topologically mixing implies uniqueness of the maximum measure for
a Q.F.T.?
 \item Is entropy a
complete invariant w.r.t. almost topological conjugacy within
topologically mixing Q.F.T.?
 \item Can one state and prove a "disjointness" property of "irreducible"
 Q.F.T. from zero-entropy systems?
 \item Does this result (or an analogue based on Yoccoz puzzle)
 extend to not necessarily connected entropy-expanding maps?
 \item Can it be applied to non-uniformly {\em hyperbolic} dynamics (by
 opposition to the non-uniformly expanding examples given here)? \footnote{We have some
 preliminary results \cite{SDPH} for partially hyperbolic diffeomorphisms with
 $\dim E^{cu}=1$.}
\end{itemize}m

\subsection*{Outline of the paper}

We first relate Q.F.T. with other classes of dynamical systems
(section \ref{sec-constructions}) before proving some basic
properties (section \ref{sec-proof-basic}). The rest of the paper
is devoted to the proof of the theorems. We first introduce the
complete Markov diagram and prove that it is conjugate with a
subset of the Q.F.T. (section \ref{sec-isom}). Then we control
measures and periodic points supported in the complement of this
subset (section \ref{sec-control}). We bound entropy at infinity
by $\hcs(\Sigma)$ (section \ref{sec-infty}). We finally prove both
theorems (section \ref{sec-theos}).

In an appendix, we prove and analyze a weaker construction
involving Hofbauer's Markov diagram instead of the complete one.

\section*{Acknowledgments}

I wish to thank Fran\c{c}ois Blanchard, Mike Boyle, Pascal Hubert
and Omri Sarig for very valuable comments/discussions. I also
thank the Ecole polytechnique f\'ed\'erale de Lausanne and the
Universit\'e de Lausanne where this work was initiated for their
warm hospitality.

\section{Construction of examples and comparisons} \label{sec-constructions}

\begin{demof}{Lemma \ref{lem-sofic}}
If $\Sigma$ is an S.F.T. then the left constraints are only the
trivial ones: the one-letter words, so that $\hcs(\Sigma)=0$.

Consider the sofic subshift over the alphabet $\{0,1,2\}$ defined
by the condition "only an even number of non-zero symbols may
appear between two $0$". Its left constraint entropy is non-zero
(it is $\log 2$).

Let $\Sigma$ be a sofic subshift with non-zero entropy. We assume
that it is irreducible (the general case follows easily). There
exists a synchronizing word $v$, i.e., such that
$\fol(uvw)=\fol(vw)$ for all words $u,w$ \cite[3.3.16]{LM}.
Therefore no left constraint can contain an occurence of $v$
anywhere except at its very beginning. But forbidding a word from
a sofic subshift strictly decreases entropy \cite[4.4.9]{LM}. This
proves that $\hcs(\Sigma)<h_\top(\Sigma)$.

\end{demof}

\begin{demof}{Lemma \ref{lem-examples-qft}}
In the two examples, we have the same situation:
\begin{enumerate}
    \item there is a natural partition $P$ defining a
symbolic dynamics $\Sigma$ by
 $$
   \Sigma := \{ A\in P^\Z : \forall n\in\Z\forall k\geq0\; <A_n\dots A_{n+k}>\ne\emptyset\};
 $$
    \item all geometric cylinders are connected because they are intervals or because they are convex;
    \item the number of $n$-cylinders which meet the boundary of the image of an
element of $P$ is bounded by $C\exp n H$ with $H<h_\top(\Sigma)$.
\end{enumerate}

Let us check condition 3:
\begin{itemize}
    \item for the case of a piecewise monotonic map, there is a finite number of
boundary points, say $N$, and only twice as much $n$-cylinders can
touch these points, hence one can take $H=0$;
    \item for the case of a piecewise entropy-expanding map this
    is part of the definition.
\end{itemize}

We conclude the proof of the Lemma by showing that
$\hcs(\Sigma)\leq H$.

Take $A_{-n}\dots A_0\in\MCs(n+1)$. By definition,
 $$
   \emptyset \ne \fol(A_{-n}\dots A_0)\subsetneq \fol(A_{-n+1}\dots A_0)
$$
This implies
 $$
  \emptyset \ne  f^{n-1}(f(A_{-n})\cap<A_{-n+1}\dots
A_0>)\subsetneq f^{n-1}<A_{-n+1}\dots A_0>
 $$
Hence, $f(A_{-n})$ meets but does not cover $<A_{-n+1}\dots A_0>$.
This last set is connected. Hence it must meet the boundary of
$f(A_{-n})$. Thus $\#\MCs(n+1)\leq\# P\times C e^{H n}$. This
concludes the proof of the Lemma.
\end{demof}

\begin{demof}{Corollary \ref{coro-entropy}}
The entropy of a piecewise monotonic map is well-known to take any
nonnegative value and the entropy of the symbolic dynamics is
equal to it.
\end{demof}

\begin{demof}{Lemma \ref{lem-hc-bad}}
We consider the subshift $\Sigma\subset\{0,a,b\}^\Z$ defined by
the following exclusions. For all $k,\ell$ distinct positive
integers, we have, using a well-known notation\footnote{for
instance, $(a|b)0^\ell(a|b)$ denotes all words $w_0\dots
w_{\ell+1}$ with $w_0,w_{\ell+1}\in\{a,b\}$ and
$w_1=\dots=w_\ell=0$.}
\begin{itemize}
  \item $(a|b)0^\ell(a|b)$ and $(a|b)0^k(a|b)$ cannot both appear in the same sequence;
  \item $(0|b)a^\ell(0|b)$ cannot appear to the right of
  $(a|b)0^\ell(a|b)$.
\end{itemize}
As we have excluded only finite words, $\Sigma$ is closed and
indeed a subshift.

It is easy to see that the invariant measures on $\Sigma$ are
supported by $\bigcup_{k\geq0} \Sigma_n'$ defined as follows.
$\Sigma_0=\{a,b\}^\Z$, $\Sigma_n$ is obtained from $\Sigma_0$ by
replacing every instance of $ba^nb$ by $b0^nb$ and $\Sigma_n'$ is
obtained by taking all sequences in $\Sigma_n$ and, for each
couple $(n,m)\in(\Z\cup\{-\infty,+\infty\})^2$, setting to $0$ all
symbols with index $\leq n$ or $>m$. Therefore, $\Sigma$ has
infinitely many maximum measures.

\medbreak

On the other hand, there is no extendable left constraint except
$O^n$, $n\geq1$, so that $\hc(\Sigma)=0$. Thus $\Sigma$ is
weak-Q.F.T. but this does not imply finiteness.
\end{demof}

\begin{demof}{Lemma \ref{lem-e-M}}
If there is an extendable left constraint, one can build an
infinite sequence $A\in\Sigma$ such that: $A_{-n}\dots A_0$ is a
left constraint for infinitely many $n$. But this means that
 $$
   \fol(A_{-n}\dots A_0)
 $$
decreases infinitely many times. Therefore we have found a
sequence in $\Sigma$ which is not eventually Markovian on the
left.
\end{demof}

\begin{demof}{Lemma \ref{lem-R}}
Fix some large integer $N\geq1$ and consider the bi-infinite
sequences obtained by concatenating blocks of the following form:
 \begin{equation}\label{eq-block1}
    (\{|[|\}|])(1|2|\dots|N)^n(a|b)^n \qquad\qquad (n\geq1)
 \end{equation}
under the constraint that matching parenthesis are of the same
type (i.e., $\{$ with $\}$, etc.). Taking the closure (which only
adds sequences of the form $(1|2|\dots|N)^\infty(a|b)^\infty$,
$(1|2|\dots|N)^\infty$ and $(a|b)^\infty$), we obtain a subshift
$\Sigma$.

The left constraints of $\Sigma$ are the blocks (we omit the
trivial, one-letter words as we shall do without further notice in
the sequel):
 \begin{itemize}
   \item $(1|2|\dots|N)^n(a|b)^k$ for $0\leq k< n<\infty$;
   \item $(\{|[|\}|])(1|2|\dots|N)^n(a|b)^k$ for $0\leq k< n<\infty$;
   \item $B_1B_2\dots B_{r-1}B_r'$ where each $B_i$ is a block
from eq. (\ref{eq-block1}), $B_r'$ is a prefix of such a block and
$B_1$ contains an opening parenthesis which is not matched in
$B_1\dots B'_r$.
 \end{itemize}
Hence $\hc(\Sigma)=\hcs(\Sigma)=\log N$.

Symmetrically, $\hc(\bar\Sigma)=\hcs(\bar\Sigma)=\frac12(\log
2+\log N)$. We see that left and right quantities are distinct.

Moreover, it is easily seen that $h_\top(\Sigma)=\log N$, $N$
being large. Thus, $\Sigma$ is a Q.F.T. with
$\hcSs(\Sigma)=\hcs(\bar\Sigma)<h_\top(\Sigma)$ but
$\hcs(\Sigma)=h_\top(\Sigma)$.
\end{demof}

\begin{demof}{Lemma \ref{lem-hc-hcs}}
The inequalities are obvious as $\MC(n)\leq\MCs(n)\leq\LL(n)$. We
describe an example where the inequalities are strict:

Take the product $\Sigma_1$ of the $2$-shift together with a
sturmian system (symbolic dynamics of a rotation by an irrational
angle $\alpha$ w.r.t. the partition
$\{[0,1-\alpha),[1-\alpha,1)\}$, see, e.g., \cite{Berthe}). Then
$\MC(\Sigma_1,n)=\MCs(\Sigma_1,n)=\LL(\Sigma_1,n)$ and
$\hc(\Sigma_1)=\hcs(\Sigma_1)=\log 2$.

\smallbreak $\Sigma_2$ will be the product of the usual even-shift
with the full $3$-shift, i.e., the subshift of
$(\{0,1\}\times\{a,b,c\})^\Z$ defined by forbidding the
words\footnote{The stars stand for any of the three symbols
$a,b,c$.}:
 $$
             (0,*)(1,*)^{2n+1}(0,*).
 $$
It is a sofic subshift hence (cf. Lemma \ref{lem-e-M})
$\hc(\Sigma_2)=0$. We compute:
 $$
   \MCs(\Sigma_2,n)=\{ (0,A_0)(1,A_1)(1,A_2)\dots(1,A_{n-1}) : A_0\dots
A_{n-1}\in\{a,b,c\}^n \}
 $$
Hence $\hcs(\Sigma_2)=\log 3<h_\top(\Sigma_2)$. Taking
$\Sigma=\Sigma_1\cup\Sigma_2$ and recalling Lemma
\ref{lem-constructions} we obtain:
 $$
    0<\hc(\Sigma)=\log 2<\hcs(\Sigma)=\log 3<h_\top(\Sigma).
 $$
\end{demof}

\begin{demof}{Lemma \ref{lem-good-coded}}
The proof is the same as the proof of Lemma \ref{lem-sofic}.
\end{demof}

\begin{demof}{Lemma \ref{lem-not-coded}}
Obvious from the remarks above the statement of the Lemma.
\end{demof}

\begin{demof}{Lemma \ref{lem-K-cap}}
We first prove $\hcs(\Sigma)\leq\text{cap}(\Sigma)$. Let $A_0\dots
A_n$ be a left constraint. Therefore one can find a finite word
$A_{n+1}\dots A_p$ such that:
 $$
  [A_0\dots A_p]=\emptyset \text{ but }
  [A_1\dots A_p]\ne\emptyset.
 $$
This implies:
 $$
   [A_1\dots A_n]_+\cap (\Sigma_+\setminus\sigma_+[A_0]_+)
   \ne\emptyset
 $$
whereas it is obvious that:
 $$
[A_1\dots A_n]_+\cap \sigma_+[A_0]_+ \ne\emptyset.
 $$
Hence, $A_1\dots A_n$ is a word that gets counted in Keller's
boundary capacity. This implies the claimed inequality.

We show that the inequality can be strict. Consider the subshift
defined by concatenating the following blocks:
 \begin{itemize}
  \item $(1|2)^n$, for any $n\geq1$;
  \item $0^nww$, for any $n\geq1$ and $w$ of the form
  $(1|2)^n$.
 \end{itemize}
The left constraints are the blocks of the form: $0^n(1|2)^k$ and
$(1|2)0^n(1|2)^k$ with $0\leq k<n$. Hence, $\hcs(\Sigma)=\log2/2$.

 Then, for all $w$ of the form $0(1|2)^n$, $n\geq0$:
 $$
    [w]_+\cap\sigma_+^n([0^{n+1}]_+)\ne\emptyset \text{ and }
    [w]_+\cap(\Sigma_+\setminus\sigma_+^n([0^{n+1}]_+)\ne\emptyset.
 $$
 so that $\text{cap}(\Sigma)=\log 2>\hcs(\Sigma)$.
\end{demof}

\begin{demof}{Lemma \ref{lem-forbid}}
If $A_0\dots A_n$ is a minimal forbidden word, then $A_0\dots
A_{n-1}$ is certainly a left constraint since it cannot be
followed by $A_n$, whereas $A_1\dots A_n$ is allowed. Thus
$\#\MCs(n)$ is at least the number of forbidden word of length
$n+1$ divided by $\#\A$. This proves that $h_{\mathcal
M}(\Sigma)\leq \hcs(\Sigma)$.

\medbreak

 We give an example where this inequality is strict.
Consider the subshift $\Sigma_1$ over $\A_1=\{0,1,2,a,b\}$ defined
by the concatenations of the following finite sequences:
 $$
        (0|1|2)^n(a|b)^k \qquad \qquad\forall n\geq1\; \forall 1\leq
k\leq n^2.
 $$
The minimal forbidden words are:
 $$
       (a|b)(0|1|2)^n(a|b)^{n^2+1} \qquad \qquad\forall n\geq1.
 $$
Hence, $h_{\mathcal M}(\Sigma_1)=\log 2$. On the other hand, the
left constraints are:
 $$
    (a|b) (0|1|2)^n(a|b)^k \qquad \qquad\forall n\geq1\; \forall 0\leq
k\leq n^2
 $$
together with the same without the first symbol $a$ or $b$.

Hence, $\hcs(\Sigma_1)=\log3$. Finally, all sequences are
eventually Markovian hence $\hc(\Sigma_1)=0$. Thus, we have:
 $$
    0 = \hc(\Sigma_1)<h_{\mathcal M}(\Sigma_1)<\hcs(\Sigma_1).
 $$

\medbreak

We now exhibit another subshift $\Sigma_2$ with
$\hc(\Sigma_2)>h_{\mathcal M}(\Sigma_2)$. $\Sigma_2$ will be
obtained by the concatenations of blocks of the same structure as
above but we introduce new, long-range restrictions to create many
extendable left constraints.

We proceed as follows. First, we restrict the blocks to $n\geq 1000$.
Then we consider:
\begin{itemize}
  \item blocks of the form $B(n):=(0|1|2)^n(a|b)^{n^2}$ with $n\geq 1000$ and even
    to be an ``opening parenthesis'' of type $[n/2]$;
  \item similar blocks but with $n$ odd to be a ``closing
    parenthesis'' of type $[n/2]$;
 \item all other blocks (i.e., all blocks with $k<n^2$) to be ``absorbing''.
\end{itemize}
The restriction is that two matching parenthesis must be of the same type
unless there is one absorbing block between them.

Thus among the left constraints are all the blocks of the form:
 $$
   (a|b)B(n_1)B(n_2)\dots B(n_r) (0|1|2)^n (a|b)^k
 $$
for all $n\geq1000$ and $1\leq k<n^2$ with $n_1,\dots,n_r$
($r\geq1$) positive integers with the restriction that $B(n_1)$ is
an opening parenthesis which is not matched and all the matchings
between $B(n_2),\dots,B(n_r)$ are between parenthesis of the same
type. It follows that $\hc(\Sigma_2)=\log3$.

On the other hand, the minimal forbidden words can be split into:
 \begin{itemize}
  \item the same as for $\Sigma_1$;
  \item $(a|b)(0|1|2)^n(a|b)$ with $n < 1000$;
  \item $B(n_1)B(n_2)\dots B(n_r)$ with $B(n_1)$ and $B(n_r)$ matching
  parenthesis of distinct types ---the point here is that only blocks with
  $n_i^2$-blocks of $a,b$ can appear.
 \end{itemize}
It follows that $h_{\mathcal M}(\Sigma_2)$ may only be slightly
larger than $\log 2$. Hence we have:
  $$
      \log 2 \approx h_{\mathcal M}(\Sigma_2) < \hc(\Sigma_2)=\hcs(\Sigma_2)=\log 3.
  $$
\end{demof}

\begin{demof}{Lemma \ref{lem-hF-bad}}
We build a subshift over $\{0,a,b\}$. Let $\X= (\X_n)_{n\in\N}$ be
a sequence of finite sets of finite words. Define
 $\Sigma_\X\subset\{0,a,b\}^\Z$ as the set of sequences such that
for all $n=1,2,\dots$, no word in $\X_n$ appears to the right of
any occurence of $(a|b)0^n(a|b)$. It is easy to check that
$\Sigma_\X$ is indeed a subshift (i.e., it is closed).

Observe that all $\sigma$-invariant probability measures of
$\Sigma_\X$ live on:
 $$
    \bigcup_{S\subset\N} \Sigma_{\X,S}
 $$
where $\Sigma_{\X,S}$ is the S.F.T. defined by excluding the words
$0^n$ for all $n\notin S$ as well as the words in $\bigcup_{n\in
S}\X_n$.

The left constraints of $\Sigma_\X$ are the (legal) words of the
following form:
 \begin{itemize}
  \item $(a|b)0^nw$ where $w$ starts with $a$ or $b$ and does not contain
  a word of the form $(a|b)0^n(a|b)$;
  \item $0^nw$ where $w$ starts with $a$ or $b$ and does not
  contain $0^n$.
 \end{itemize}

We set:
 $$
   \X_n=\{(0|b)a^n(0|b)\} \qquad\forall n\geq1.
 $$
We have $\Sigma_{\X,\emptyset}=\{a,b\}^\Z$ and each
$\Sigma_{\X,S}\subset\{0,a,b\}^\Z$ is obtained from $\{a,b\}^\Z$
by substituting $b0^nb$ for all blocks $ba^nb$, for $n\in N$.
Hence for all $S\subset\N$
 $$
    h_\top(\Sigma_{\X,S})= \log 2 = h_\top(\Sigma).
 $$
We see that there are {\bf uncountably many maximum measures },
one on each S.F.T. $\Sigma_{\X,S}$, $S\subset\N$.

We have $\hc(\Sigma)=\hcs(\Sigma)=h_\top(\Sigma)$. Indeed, it
follows from the main theorem (or can be easily checked from the
above description of left constraints and the observation that any
minimal left constraint can be extended by inserting longer and
longer runs of $0$, so that $\hc(\Sigma)=\hcs(\Sigma)$).

We claim that $h_\fol(\Sigma)=0$. The follower set of a word $w$
is described by giving the set of distinct lengths of the
$0$-blocks bounded by letters in $w$ (say
$1\leq\ell_1<\ell_2<\dots<\ell_r$ for some $0\leq r<n$) together
with the lengths $0\leq \ell_+,\ell_-\leq n$ of the runs of zeroes
that begin and end $w$.

We see that $n\geq \sum_{i=1}^r\ell_i\geq\sum_{i=1}^r i\geq
r^2/2$. Hence $r\leq2\sqrt{n}$. Therefore the number of distinct
follower sets defined by words of length $n$ is bounded by:
$(n+1)^2 C^{2\sqrt{n}}_n$. This proves the claim.
\end{demof}

\section{Proofs of basic properties} \label{sec-proof-basic}

\begin{demof}{Proposition \ref{prop-top-invariance}}
We first observe that the Q.F.T. and weak-Q.F.T. properties are
not preserved under extensions or factor maps, already for trivial
reasons:

Indeed, take a Q.F.T. $\Sigma_1$ and a subshift $\Sigma_2$ with
the same entropy which is not a Q.F.T. Let $\pi|\Sigma_1=\id$
whereas $\pi(\Sigma_2)$ is a fixed point. Then consider
$\pi:\Sigma_1\cup\Sigma_2\to\Sigma_1\cup\{0\}$:
$\Sigma_1\cup\{0\}$ is a Q.F.T. with an extension,
$\Sigma_1\cup\Sigma_2$ which is not.

Take now a Q.F.T. $\Sigma_1$ and a subshift $\Sigma_2$ with a
strictly smaller entropy which is not a Q.F.T. Let
$\pi|\Sigma_2=\id$ whereas $\pi(\Sigma_1)$ is a fixed point $0$.
Then consider $\pi:\Sigma_1\cup\Sigma_2\to\Sigma_2\cup\{0\}$:
$\Sigma_1\cup\Sigma_2$ is a Q.F.T. with an image
$\{0\}\cup\Sigma_2$ which is not.

Now, we have seen that there are sofic subshifts $\Sigma$ with
$\hcs(\Sigma)>0$ whereas of course their S.F.T. extension
$\Sigma_0$ has $\hcs(\Sigma_0)=0$ so that the left constraint
entropy does not always decrease under factor maps.

We finally turn to the invariance of $\hcs(\Sigma)$. Let
$h:\Sigma'\to\Sigma$ be the conjugacy. We have:
 $$
    (h(x))_i = H(x_{i-L}\dots x_{i+L}) \qquad\forall
    x\in\Sigma'\;\forall i\in\Z
 $$
for some integer $L\geq0$ and some map $H:{\A'}^{2L+1}\to\A$.
Similarly, there is a map $H':\A^{2L+1}\to\A'$ for $h^{-1}$ (maybe
after increasing $L$).

For all $n$ large enough, we shall construct a map $\psi:
\MCs(\Sigma,n)\to\bigcup_{k=-L}^L \MCs(\Sigma',n+k)$ which is at
most $\#{\A'}^{2L}$ to $1$. This will clearly imply
$\hcs(\Sigma)\leq \hcs(\Sigma')$ (notice that this would not work
if $\hcs(\Sigma)$ were defined using $\liminf$ instead of
$\limsup$).

Thus we take $A_{-n}\dots A_0$ a left constraint of $\Sigma$. We observe that
there exist $L,R\in\Sigma$ such that for some $0\leq k<\infty$:
 \begin{equation}\label{eq-lc}
    L_{-n}\dots L_0 = A_{-n}\dots A_0,\quad
    R_{-n+1}\dots R_0 = A_{-n+1}\dots A_0,\text{ and }
    [A_{-n}R_{-n+1}\dots R_k]=\emptyset.
\end{equation}

Let $L'=h^{-1}(L)$ and $R'=h^{-1}(R)$. Observe that
 $$
   L'_{-n+L+1}\dots L'_{-L} = R'_{-n+L+1}\dots R'_{-L}.
 $$

\begin{claim}
There is $0\leq\ell\leq 2L$ such that $L'_{-n+L-\ell} \dots
L'_{-L}$ is a left constraint.
\end{claim}

The claim will give the map $\psi$ discussed above and therefore
the inequality for $\hcs(\Sigma)$.

We prove the claim by contradiction. We first observe that
$\fol(L'_{-n+L-\ell}\dots L'_{-L})$ is evidently non-empty for
all $\ell\geq0$. Hence, if the claim is false, it means
that:
 $$
    \fol(L'_{-n-L} \dots L'_{-L})
       =\fol(L'_{-n+L+1}\dots L'_{-L})
       =\fol(R'_{-n+L+1}\dots R'_{-L})\ni R'_{-L}R'_{-L+1}\dots
 $$
Hence $[L'_{-n-L} \dots L'_{-n+L}R'_{-n+L+1}\dots R'_{k+L}]\ne
\emptyset$. Applying $h$, we find that $[L_{-n}R_{-n+1}\dots
R_k]\ne\emptyset$, a contradiction. The claim is proved.

\medbreak

We now turn to $\hc(\Sigma)$.

Let $A_{-n}\dots A_0$, $n\geq0$, be an extendable left constraint.
We can find $L\in\Sigma$ with $L_{-n}\dots L_0=A_{-n}\dots A_0$
such that for infinitely many integers $m\geq0$, there exist
$R^{(m)}\in\Sigma$ and $k^{(m)}\geq0$ such that:
 $$
    R^{(m)}_{-m+1}\dots R^{(m)}_0 = L_{-m+1}\dots L_0 \text{ and }
      [L_{-m}\dots L_0 R^{(m)}_1\dots
      R^{(m)}_{k^{(m)}}]=\emptyset.
 $$
 Applying the previous argument
we obtain for each value of $m$, a left constraint
$L'_{-m+L-\ell(m)} \dots L'_{-n+L}$ with $0\leq\ell(m)\leq 2L$.
Hence we see that $L'_{-m+L-\ell} \dots L'_{-n+L}$ is indeed an
extendable left constraint.
\end{demof}

\begin{demof}{Lemma \ref{lem-constructions}}
Let $X$, resp. $Y$, be a subshift over the alphabet $\A$, resp.
$\mathcal B$. We claim that
 $$
  \fol(A_{-n},B_{-n})\dots (A_0,B_0)) = \fol(A_{-n}\dots A_0)
   \times \fol(B_{-n}\dots B_0).
 $$
Indeed, observe that:
 \begin{equation}\label{eq-prod}
   (A_{-n},B_{-n})\dots (A_0,B_0) \in\MCs(X\times Y,n) \iff
  A_{-n}\dots A_0\in\MCs(X,n) \text{ or }
  B_{-n}\dots B_0\in\MCs(Y,n).
 \end{equation}
Thus,
 $$
    \MCs(X\times Y,n) = \MCs(X,n)\times\LL(Y,n) \cup
         \LL(X,n)\times \MCs(Y,n)
 $$
so that $\hcs(X\times Y) = \max(\hcs(X)+h_\top(Y),
h_\top(X)+\hcs(Y))$. This gives the result for $\hcs$.

But it is obvious that the equivalence (\ref{eq-prod}) is also
valid for extendable left constraints. This concludes the proof of
the Lemma.
\end{demof}

\section{Partial conjugacy} \label{sec-isom}

We shall build a conjugacy with the following system:

\begin{defi}
The \new{complete Markov diagram} of $\Sigma$ is the graph $\MD$
the vertices of which are the left constraints and the arrows:
$A_{-n}\dots A_0\to B_{-m}\dots B_0B_1$ if and only if: $m\leq n$
and
 $$
  B_{-m}\dots B_0=A_{-m}\dots A_0 \text{ and } \fol(A_{-n}\dots A_0B_1)= \fol(B_{-m}\dots B_0)
 $$
The corresponding Markov shift is denoted by $\Sh$.
\end{defi}

\begin{remk}
This is a variant of Hofbauer's Markov diagram. However, it is
necessary to use this variant to exploit the bound on
$\hc(\Sigma)$. See the Appendix.
\end{remk}

\subsection*{Partial isomorphism}

The natural projection $\pi:\Sh\to\Sigma$ is defined by
$(\pi(\alpha))_n=A$ iff the finite word $\alpha_n$ ends in $A$.

\begin{lem}
$\pi:\Sh\to\Sigma$ is well-defined.
\end{lem}

\begin{demo}
Let $\alpha\in\Sh$ and set $A=\pi(\alpha)$. We have to prove that
for all $n\in\Z$, $p\geq0$,
 $$
    [A_n\dots A_{n+p}]\ne \emptyset.
 $$
But it follows from the definition of the arrows of $\MD$ and an
immediate induction that:
 \begin{equation} \label{eq-alpha}
   \fol(\alpha_{n+p})=\sigma_+^p(\fol(\alpha_n)\cap[A_n\dots A_{n+p}])
 \end{equation}
As $\fol(\alpha_{n+p})\ne\emptyset$, the lemma is proved.
\end{demo}

 \medbreak

The conjugacy will be restricted to a set $\Sigma_M\subset\Sigma$.
Recall Definition \ref{def-e-M} of an eventually Markovian.

\begin{defi}
$A\in\Sigma$ is \new{completely Markovian} iff $\sigma^nA$ is
eventually Markovian for all $n\in\Z$.

The set of completely Markovian sequences is denoted by
 $\Sigma_M$.
\end{defi}

\begin{prop} \label{prop-conjugacy-complete}
The restriction $\pi:\Sh\to\Sigma_M$ is a conjugacy.
\end{prop}

\begin{demo}
We define a partial inverse $i:\Sigma_M\to\Sh$ to $\pi$ by the
formula:
 $$
   i(A)=\alpha \text{ with } \alpha_n = A_{n-\ell}\dots A_n
 $$
where $\ell=\ell(A,n)$ is the minimum integer such that, for all
$k\geq\ell$,
 $$
   \fol(A_{n-k}\dots A_n)=\fol(A_{n-\ell}\dots A_n).
 $$
We check that for all $A\in\Sigma_M$ $i(A)$ is a well-defined
element of $\Sh$:
\begin{itemize}
 \item As $\ell$ is chosen minimum, $A_{n-\ell}\dots A_n$ is
indeed a left constraint, hence a vertex of $\MD$;
 \item Taking
$L\geq\max(\ell(n),\ell(n-1)+1)$, we have:
  $$
   \fol(A_{n-\ell(n)}\dots A_n) = \fol( A_{n-L}\dots
A_{n-1}A_n)
  $$
  $$
  \qquad\qquad\qquad = \sigma_+(\fol(A_{n-L}\dots A_{n-1}))\cap A_n
  $$
  $$
  \qquad\qquad\qquad = \fol( A_{n-1-\ell(n-1)}\dots
A_{n-1}A_n),
 $$
hence $(i(A))_{n-1}\to (i(A))_n$ is an arrow of $\MD$.
\end{itemize}

It is clear that $\pi\circ i = \id_{\Sigma_M}$.

It remains to see that $\pi(\Sh)\subset\Sigma_M$. Let
$\alpha\in\Sh$ and $A=\pi(\alpha)$. We have to prove that
$\sigma^nA$ is eventually Markovian for all $n\in\Z$. We consider
the case $n=0$, the general case being exactly the same.

Let $m$ be the length of the left constraint $\alpha_0$.
$\alpha_{-m}$ is some left constraint $C_{-p}\dots C_{-1}A_{-m}$
for some $p$. We prove by induction that:
 $$
   \alpha_{-m+k} \text{ is a suffix of } C_{-p}\dots
C_{-1}A_{-m}\dots A_{-m+k}.
 $$
Indeed it is true for $k=0$ and the definition of $\MD$ ensures
that $\alpha_{-m+k+1}$ is a suffix of $\alpha_{-m+k}A_{-m+k+1}$.

Therefore, $\alpha_0$ is the suffix of length $m$ of
 $$
C_{-p}\dots C_{-1}A_{-m}\dots A_0.
 $$
Hence $\alpha_0=A_{-m}\dots A_0$. By the same token,
 $$
   \fol(A_{-q}\dots A_0) = \fol (A_{-m}\dots A_0)
 $$
for all $q\geq m$. This proves that $A$ is eventually Markovian
and concludes the proof of the proposition.
\end{demo}

\section{Control of the non-Markovian part} \label{sec-control}

We prove the simpler statement for periodic points first:

\begin{lem} \label{lem-per-non-M}
The periodic orbits in $\Sigma\setminus\Sigma_M$ satisfy:
 $$
    \limsup_{n\to\infty}
    \frac1n\log\#\{x\in\Sigma\setminus\Sigma_M: \sigma^nx=x\} \leq
    \hc(\Sigma).
 $$
\end{lem}

\begin{demo}
Let $X_n:=\{A\in\Sigma\setminus\Sigma_M: \sigma^nA=A\}$ and take
$A\in X_n$. As $A$ is not eventually Markovian, we have, that for
infinitely many $k\geq0$, $A_{-k}\dots A_0\in\MC(\Sigma,k+1)$.
Hence, for such a $k$, $A_{-k}\dots A_{-k+n-1}\in\MC(\Sigma,n)$.
$A$ being $n$-periodic, this means that $\#X_n\leq
n\times\#\MC(\Sigma,n)$.
\end{demo}

We turn to the measures:

\begin{prop} \label{prop-non-M}
Let $\mu$ be a $\sigma$-invariant probability
measure with $\mu(\Sigma\setminus\Sigma_M)=1$. Then
 $$
          h(\mu,\sigma) \leq \hc(\Sigma).
 $$
\end{prop}

\begin{remk}
The above estimate is sharp in that the inequality can be an
equality: take the union of $3$-shift and of the product of the
$2$-shift with the symbolic dynamics of an irrational rotation.
Then there is a measure on $\Sigma_M$ with entropy $\log
2=\hc(\Sigma)<h_\top(\Sigma)=\log 3$.
\end{remk}

\begin{demo}
We fix $\mu$ as above and bound its entropy. We denote by $N$ the
set of sequences which are not eventually Markovian. We first
claim that $\mu(N)=1$. Indeed, if $A$ is not eventually Markovian,
then
 \begin{equation} \label{eq-fol-dec}
     \fol(A_{-n}\dots A_0) \subsetneq \fol (A_{-n+1}\dots A_0)
 \end{equation}
for infinitely many $n\geq0$. But eq. (\ref{eq-fol-dec}) is
equivalent to:
 $$
    \sigma_+(A_{-n})\not\supset [A_{-n+1}\dots A_0]_+
 $$
This last condition obviously implies:
 $$
    \sigma_+(A_{-n})\not\supset [A_{-n+1}\dots A_{-1}]_+.
 $$
Thus, $\sigma^{-1} A$ is also not eventually Markovian. Therefore
$\sigma^{-1}(N)\subset N$. We have
$\mu(\sigma^{-1}(N))=\mu(N)$ by the $\sigma$-invariance of $\mu$.
Thus, $\mu(N\Delta\sigma^{-1}(N))=0$. We conclude that
 $$
  \Sigma\setminus \Sigma_M=\bigcup_{n\in\Z}\sigma^n(N)
   = N
 $$
up to $\mu$-negligible sets, hence by ergodicity, $\mu(N)=1$ as
claimed. This argument is due to Hofbauer.

\medbreak

We bound the entropy of $\mu$ by bounding the minimal number of
$n$-cylinders whose union has measure $>1/2$ (see, e.g.,
\cite{DR}). Let $\eps>0$.

Let $K_0<\infty$ be such that $\#\MC(n)\leq
e^{(\hc(\Sigma)+\eps)n}$ for all $n\geq K_0$. We also assume $C^{2n/K_0}_n
\leq e^{\eps n}$ for all large $n$.

Let $n(A)=\min\{k\geq K_0:A_{-k}\dots A_0 \in\MC(k+1)\}$. As
$\mu(N)=1$, $n(\cdot)<\infty$ $\mu$-a.e.

There exists $N_0<\infty$ such that $n(\cdot)>N_0$ on a set of
measure $<\eps/\log\#\A$.

By Birkhoff's ergodic theorem, there exist an integer $M_0<\infty$
and a measurable set $G_0\subset\Sigma$ with $\mu(G_0)>1/2$ such
that for all $A\in G_0$, all $n\geq M_0$,
 $$
   \frac1n\#\{ 0\leq k<n : n(\sigma^kA)>N_0\} <
\frac{\eps}{\log\#\A}.
 $$
We may and do assume that $M_0\geq N_0 \log\#\A/\eps$.

It is easy to see that for any $n\geq M_0$, any $A\in G_0$,
$A_0\dots A_{n-1}$ can be decomposed into:
 \begin{itemize}
  \item segments belonging to some $\MC(\ell)$ with $\ell\geq
K_0$;
  \item an initial segment of length at most $N_0$;
  \item at most $\eps n/\log\#\A$ left-overs.
 \end{itemize}

Thus, the number of $n$-cylinders meeting $G_0$ is bounded by:
 $$
    C^{2n/K_0}_n e^{(\hc(\Sigma)+\eps)n} \#\A^{N_0+\eps
n/\log\#\A}
       \leq e^{(\hc(\Sigma)+4\eps)n}
 $$
for all large $n$.

As $\eps>0$ was arbitrary, this proves that $h(\mu,\sigma)\leq
\hc(\Sigma)$.
\end{demo}

\section{Entropy at infinity} \label{sec-infty}

\begin{prop} \label{prop-ent-inf}
For any $\eps>0$, there exist a number $\delta>0$ and a finite
subset $\MD_0\subset \MD$ such that any ergodic,
$\sigma$-invariant probability measure $\muh$ on $\Sh$ such that
$\muh\left(\bigcup_{D\in\MD_0}[D]\right)<\delta$ satisfies:
$h(\muh,\sigma)\leq \hcs(\Sigma)+\eps$.
\end{prop}

\begin{demo}
Let $\eps>0$. Let $K_0<\infty$ be such that for all $n\geq K_0$,
$\#\MCs(n)\leq e^{(\hcs(\Sigma)+\eps)n}$. We assume $K_0$ to be
large enough so that $C^{2n/K_0}_n\leq e^{\eps n}$ for all large
$n$.

Let $\MD_0:=\bigcup_{n\leq K_0}\MCs(n)$. Let
$0<\delta<\eps/(K_0\log\#A)$ be such that $C^{2\delta n}_n\leq
e^{\eps n}$ for all large $n$.

Let $\muh$ be as above. We bound its entropy as in the proof
of Proposition \ref{prop-non-M} by finding
an upper bound for the number of $n$-cylinders $A_0\dots A_{n-1}$
of the form $A=\pi(\alpha)$ with:
 $$
    \frac1n \#\{ 0\leq k<n: \ell(\alpha_k) < K_0 \} \leq \delta
 $$
and $\ell(\alpha_0)<L_0$ for some large $L_0$, possibly depending
on $\muh$.

We cut $A_0\dots A_{n-1}$ into maximal segments according to
whether $\ell(\alpha_k)<K_0$ or not. There are at most $2\delta n$
cutting points. Hence at most $C^{2\delta n}_n\leq e^{\eps n}$
choices of positions.

Each interval below level $K_0$ is described by giving directly
the symbols involved. There are at most $\#\A^{\delta n}\leq
e^{\eps n}$ choices.

Each interval $A_{m}\dots A_{m+k-1}$ above level $K_0$ is in turn
divided into sub-segments as follows. We start from the end
setting $n_0=m+k$ and, inductively,
$n_{i+1}=n_i-\ell(\alpha_{n_i})$. We stop at the smallest $i=i_*$
such that $n_i\leq m$. Thus, there are at most
$e^{(\hcs(\Sigma)+\eps)(m+k-n_{i_*})}$ choices of symbols.

We have to find a lower bound for $n_{i_*}$. Observe that
$\ell(\alpha_{n_{i_*-1}}) \leq \ell(\alpha_m)+n_{i_*-1}-m$.

If $m>0$, then $\ell(\alpha_m)=K_0$ and $n_{i_*}\geq m-K_0$. The
number of choices of symbols for the interval $A_{m+1}\dots
A_{m+k}$ is bounded by $e^{(\hcs(\Sigma)+\eps)(k+K_0)}$.

If $m=0$, then $\ell(\alpha_0)\leq L_0$ and the number of choices
of symbol is bounded by $e^{(\hcs(\Sigma)+\eps)(k+L_0)}$.

 We notice that there are at most $\delta n+1$ such intervals.

Taking product, we find a total number choices for $A_0\dots
A_{n-1}$ bounded by
 $$
   e^{\eps n} e^{(\hc(\Sigma)+\eps)L_0} e^{\hcs(\Sigma)+\eps)n}
    e^{(\hcs(\Sigma)+\eps)K_0(\delta n+ 1)} \leq
C e^{2\eps n} e^{(\hcs(\Sigma)+\eps)(1+\delta)n}
 $$
with $C= e^{(\hcs(\Sigma)+\eps)(L_0+K_0)}$. Thus,
 $$
   h(\muh,\sigma)\leq (\hcs(\Sigma)+4\eps).
 $$
But $\eps>0$ was arbitrary.
\end{demo}

\section{Proof of the Theorems} \label{sec-theos}

We may and do assume that $\hcs(\Sigma)=\hcSs(\Sigma)$ (or
$\hc(\Sigma)=\hcS(\Sigma)$ depending on the case). Otherwise
replace $\Sigma$ by $\bar{\Sigma}$.

\subsection{Structure theorem}

We collect the previous results that imply the structure theorem.

The countable oriented graph $\mathcal G$ of the statement is the
complete Markov diagram. $\pi:\Sigma(\mathcal G)\to\Sigma$ is induced by
the natural projection $\mathcal G\to\mathcal A$.
$\mathcal X$ is $\Sigma\setminus\Sigma_M$. Then the conjugacy
between $\Sigma(\mathcal G)$ and $\Sigma\setminus\mathcal X$ is
given by Proposition \ref{prop-conjugacy-complete}. The control on
$\mathcal X$ follows from Proposition \ref{prop-non-M} (for
measures) and Lemma \ref{lem-per-non-M} (for periodic points). The
finiteness at infinity follows from Proposition
\ref{prop-ent-inf}.

This concludes the proof of the structure theorem.

\subsection{Main theorem}

We deduce Theorem \ref{theo-main} from the structure theorem,
Theorem \ref{theo-struc}. To begin with, we consider the case of a
weak-Q.F.T.

The first point follows from the conjugacy (up to measures of
entropy $\leq \hcS(\Sigma)$) using the following result of
Gurevich \cite{Gurevich}: on each irreducible Markov shift there
is at most one maximum measure --and obviously there are at most
countably irreducible Markov subshifts, as the graph itself is
countable. The maximum measures are Bernoulli by Proposition 2 of
section 5 of \cite{IMT} which shows that Markovian measures are
weak Bernoulli and therefore Bernoulli by Ornstein's isomorphism
theorem.

Let us prove the estimate on the number of periodic points. The
upper bound follows from the definition of topological entropy. We
establish the lower bound. First observe that $\Sigma$, as a subshift
on a finite alphabet is expansive and therefore admits a maximum measure
(see \cite{DGS}). 

Using the conjugacy of the structure theorem and the
variational principle of Gurevich  (see \cite{Kitchens} for background on
Markov shifts), we see that 
$\mathcal G$ has an irreducible subgraph defining a Markov shift which carries
a probability measure with entropy equal to its Gurevich 
entropy equal to $h_\top(\Sigma)$.
This implies that the the number $\ell_n$ of loops of length $n$ at a given
vertex in this irreducible subgraph satisfies
$\limsup_{n\to\infty} \frac1n\log\ell_n= h_\top(\Sigma)$ according to
Vere-Jones. 

But these loops define periodic points for $\Sigma(\mathcal G)$ hence for 
$\Sigma$ using the embedding $\pi$. This proves the estimate on periodic 
points.

\medbreak

We now turn to the case of a Q.F.T.

The finite number of maximum measures follows from Proposition
\ref{prop-ent-inf}.

Deferring the proof on the meromorphy of the zeta function we
recall how to deduce the estimate on periodic points from it
(one usually finds more delicate estimates, see, e.g., page 101
of \cite{PP}).

There are $1\leq t\leq s<\infty$, $\kappa<1$, complex numbers
$\lambda_1,\dots,\lambda_s$ with moduli $e^{h_\top(\Sigma)}$, and
positive integers $q_1,\dots,q_s$ such that:
 $$
     \zeta(z) = \psi(z) \prod_{i=1}^s (1-\lambda_i^{-1}z)^{-q_i}
 $$
where $\psi(z)$ is holomorphic and non-zero on
$|z|<\rho$ for some $\rho>e^{-h_\top(\Sigma)}/\kappa$. 
We compute the logarithmic derivative of each side:
 $$
    \frac{\zeta'(z)}{\zeta(z)} =
      -\sum_{n\geq0} \#\{x\in\Sigma:\sigma^{n+1} x=x\}z^{n}
 $$
and
 $$
    \frac{\psi'(z)}{\psi(z)}
      - \sum_{i=1}^s \frac{\lambda_i^{-1}q_i}{1-\lambda_i^{-1}z}
    = \sum_{n\geq0} \left(\phi_n-\sum_{i=1}^s
    \lambda_i^{-n-1}q_i\right) z^n
 $$
with $|\phi_n|\leq C\kappa^n e^{n h_\top(\Sigma)}$ as
$\psi'(z)/\psi(z)$ is analytic on
$|z|<\rho$. It follows that:
 $$
    \left|e^{-(n-1)h_\top(\Sigma)}\#\{x\in\Sigma:\sigma^n x=x\} - \sum_{i=1}^s
    e^{n\sqrt{-1}\theta_i}q_i \right|
      \leq C\kappa^n
 $$
where $e^{\sqrt{-1}\theta_i}=\lambda_i/|\lambda_i|$. But, as $q_i\ne0$, it
follows that
 $$
   0<\limsup_{n\to\infty} \sum_{i=1}^s
    e^{\sqrt{-1}n\theta_i}q_i < \infty
 $$
The claim on the number of periodic points follows.

\medbreak

 It remains to prove the analyticity properties of the
zeta function.

\subsection*{The zeta function}

Lemma \ref{lem-per-non-M} immediately implies that
$\frac{\zeta_\Sigma(z)}{\zeta_{\Sigma_M}(z)}$ is a holomorphic
function over the disk $|z|<e^{-\hc(\Sigma)}$. Hence it is enough
to prove the result for $\zeta_{\Sigma_M}(z)$.

Observe that, $\pi$ being a conjugacy between $\Sh$ and
$\Sigma_M$, it defines a bijection between periodic points of
$\Sh$ and of $\Sigma_M$ and this bijection of course preserves
minimum period. Hence $\zeta_{\Sigma_M}(z)=\zeta_{\Sh}(z)$. We
study this last function adapting the proof of Hofbauer and Keller
from \cite{HK}.

\medbreak

The Markov diagram $\MD$ defines a countable matrix
$K:\MD\times\MD\to \{0,1\}$ according to: $K(i,j)=1\iff i\to j$.
We observe that:
 $$
    \zeta_{\Sh}(z) = \exp -\sum_{n\geq1} \frac{z^n}n \tr K^n
 $$
where $\tr K^n:=\sum_{i\in\MD} (K^n)(i,i)$. We observe that for
each $n$, $\tr K^n<\infty$. In fact it is bounded by
$\#\LL(\Sigma,n)$ which is at most of the order of $e^{n
h_\top(\Sigma)}$. This proves the analyticity claim.

For convenience we assume some identification of $\MD$ with $\N$.
Let $n\geq1$ be some integer. Write $K=\left({ A \; U \atop V \;
B}\right)$ where $A$ is a $n\times n$-submatrix.

Let $k>n$ be some other integer. Write $\tilde{K}$ for the {\bf
finite} matrix obtained by truncation of $K$ to the indices
$(i,j)$ with $\max(i,j)\leq k$. Write $\tilde{K}=\left({ A \;
\tilde{U} \atop \tilde{V} \; \tilde{B}}\right)$.

\begin{claim}{\bf A.}
Given $\eps>0$, the spectral radius of $\tilde{B}$,
$\rho(\tilde{B})$, is bounded by $e^{\hcs(\Sigma)+\eps}$ as soon
as $n\geq n_0(\eps)$.
\end{claim}

Indeed, each coefficient $(\tilde{B}^m)_{ij}$ is bounded by the
number of paths of length $m$ on the subset of $\MD$ corresponding
to the integers $\geq n$ and starting at $i$ and ending at $j$.
Therefore (cf. the proof of Proposition \ref{prop-ent-inf})
$\|\tilde{B}^m\|$ grows at most like
$C_{n_0}(i)e^{(\hc(\Sigma)+\eps)m}$ if $n$ is large enough. This
proves the claim.

\medbreak

\begin{lem} \cite[Lemma 2]{HK}
If $L$ is a finite matrix and if $L=\left({ L_{11} \; L_{12} \atop
L_{21} \; L_{22}}\right)$ is a block decomposition with $L_{22}$
invertible, then:
 $$
     \det L = \det L_{22} \det(L_{11}-L_{12}L_{22}^{-1}L_{21})
 $$
\end{lem}

We apply this Lemma to $I-z\tilde{K}=\left({ I-zA \;\; -z\tilde{U}
\atop -z\tilde{V} \;\; I-z\tilde{B}}\right)$ ($I$ denotes each
time the identity matrix of the required dimensions). This is
possible because, by the previous claim, for all
$|z|<e^{-\hcs(\Sigma)}$, $I-z\tilde{B}$ is invertible. Thus,
 $$
    \det(I-z\tilde{K}) =
      \exp \det(I-z\tilde{B})
      \times
      \det(I-zA-z^2\tilde{U}(I-z\tilde{B})^{-1}\tilde{V}).
 $$
 Hence
  $$
 \exp \left(-\sum_{m\geq1} \frac{z^m}{m} \tr\tilde{K}^m\right) =
      \exp \left(-\sum_{m\geq1} \frac{z^m}{m} \tr\tilde{B}^m\right)
      \times
      \det(I-zA-z^2\tilde{U}(I-z\tilde{B})^{-1}\tilde{V}).
 $$

\begin{claim}{\bf B.}
We have uniform convergence on all compact subsets of
$|z|<e^{-h_\top(\Sigma)}$ when $k\to\infty$, of
 $$
  \exp \left(\sum_{m\geq1} \frac{z^m}{m} \tr\tilde{K}^m\right)
     \to \exp \left(\sum_{m\geq1} \frac{z^m}{m} \tr K^m\right) = \zeta_\Sh(z)^{-1}.
 $$
The same is true for $|z|<e^{-\hcs(\Sigma)}$ for $B$ instead of
$K$. Call $B_n(z)$ the resulting analytic function.
\end{claim}

The claim follows from routine arguments if
 $$
   0\leq\tr\tilde{K}^m\leq \tr K^m\leq C_\eps
   e^{(h_\top(\Sigma)+\eps)m}
 $$
and, (for $m\geq m(\eps)$)
 \begin{equation}\label{boundTrB}
   0\leq\tr\tilde{B}^m\leq \tr B^m\leq C_\eps
   e^{(\hcs(\Sigma)+\eps)m}
 \end{equation}
with constants $C_\eps$ independent of $k$.

These inequalities requires a little care since they {\it a priori} involve 
infinitely many coefficients, in contrast to Claim A.

\medbreak

For the first inequality, it is enough to remark that
 $$
   \tr K^m \leq \#\{x\in\Sigma:\sigma^mx=x\}
 $$
using that $\pi$ embeds $\Sh$ into $\Sigma$. The inequality follows from the
definition of the entropy of $\Sigma$.

\medbreak

We turn to $\tr B^m$. It is obviously bounded by the number of closed paths of
length $m$ which stay above $n$ in $\MD$. Observe also that, by Proposition
\ref{prop-conjugacy-complete}, it is enough to count the
projections on $\Sigma$ of these loops.

Take one such loop. It determines $\alpha\in\Sh$ with
$\sigma^m\alpha=\alpha$. $\alpha$ projects to $A\in\Sigma$ with
$\sigma^mA=A$. Set $i_0=m$ and define recursively
$i_{s+1}=i_s-\ell(\alpha_{i_s-1})$ for $s\geq0$ ($\ell(\gamma)$ is the length 
of the finite word $\gamma$) . Let $S$ be the
smallest integer such that $i_S\leq0$. We have cut the sequence
$A_{i_S}\dots A_{m-1}$ into left constraints with length $\geq L(n):=
\min\{\ell(\gamma):\gamma\in\MD\text{ with }\gamma\geq n\}$.

Consider first the case where $i_S\leq -m$. Recall that any prefix
of a left constraint is a left constraint. Hence $A_{i_S}\dots
A_{i_S+m-1}\in\MCs(m)$. Such loops of length $m$ are therefore in
numbers bounded by $m C_\eps e^{m (\hcs(\Sigma)+\eps)}$, as
pbviously $\lim_{n\to\infty} L(n)=\infty$.

Now assume that $-m<i_S\leq0$, so that $0<i_S+m\leq m$. Set $t:=1$ if $i_S=0$
or $t:=\min\{s:i_s<i_S+m\}$ otherwise. We consider the
cutting of $A_{i_S}\dots A_{i_S+m-1}$ into
 $$
 A_{i_S}\dots A_{i_S-1},
\dots,A_{i_{t-1}}\dots A_{i_t-1},A_{i_t}\dots A_{i_S+m-1}
 $$
where each block is a left constraint and each has a length at
least $L(n)$, except possibly the last. It is now easy to bound the
number of such loops by $C_\eps e^{m(\hcs(\Sigma)+\eps)}$ (cf. the
proof of Proposition \ref{prop-ent-inf}). This concludes the proof
of eq. (\ref{boundTrB}).

\begin{claim} {\bf C.}
We have uniform convergence on all compact subsets of
$|z|<e^{-(\hcs(\Sigma)+\eps)}$ when $k\to\infty$ of the $n\times
n$-matrices
 $$
    I-zA-z^2\tilde{U}(I-z\tilde{B})^{-1}\tilde{V} \to
     \underbrace{I-zA-z^2U(I-zB)^{-1}V}_{=:D_n(z)}.
 $$
\end{claim}

Indeed,
 $$
   \tilde{U}(I-z\tilde{B})^{-1}\tilde{V} = \sum_{m\geq0}
   z^m\tilde{U}\tilde{B}^m\tilde{V}
 $$
and the coefficient $(\tilde{U}\tilde{B}^m\tilde{V})_{i,j}$ is the
number of paths of length $m+2$ going from $i$ to $j$ with
$i,j\leq n$ staying above $n$ and below $k$. As $\MD$ has finite
outdegree, a path of length $m$ starting from $i\leq n$ cannot go
above some integer $k_1(n,m)$. Hence,
$\tilde{U}\tilde{B}^m\tilde{V}=UB^mV$ as soon as $k\geq k_1(n,m)$.

Moreover, once again for the same reasons as in the proof of Claim
A,
 $$
   0\leq(\tilde{U}\tilde{B}^m\tilde{V})_{i,j}\leq
                 C_{\eps, n_0} e^{m(\hcs(\Sigma)+\eps)}
 $$
The claim C follows immediately.

\medbreak

To conclude, we see by Claims B and C that
$\zeta_\Sh(z)^{-1}=B_n(z) \det D_n(z)$ on
$|z|<e^{-h_\top(\Sigma)}$. But the right hand side has an obvious
holomorphic extension to $|z|<e^{-(\hcs(\Sigma)+\eps)}$ by Claim C.

Finally, by letting $n\to\infty$, we obtain the result on the full
disk $|z|<e^{-\hcs(\Sigma)}$. This proves the claimed properties
of the zeta function and concludes the proof of the Theorem
\ref{theo-main}.

\section*{Appendix: Hofbauer's Markov shift}
\label{sec-Hof-diag}

\begin{defi}
The \new{Hofbauer Markov diagram} of $\Sigma$ is the graph $\MDs$
the vertices of which are the follower sets and the arrows:
 $$
   F\to G \iff \exists A\in\A\; G=\sigma_+(F)\cap[A].
 $$
The corresponding Markov shift is denoted by $\Shs$.
\end{defi}

The natural projection $\pi:\Shs\to\Sigma$ is defined by
$(\pi(\alpha))_n=A$ iff $\alpha_n\subset [A]$. $\pi$ does not
define an isomorphism of the whole of $\Shs$, but we have to take
a subset:

\begin{defi}
 $\alpha\in\Shs$ is \new{explicitely Markovian} if
for all $n\in\Z$, there exists $\ell=\ell(n)\geq0$ such that
 $$
    \alpha_n = \fol(A_{n-\ell}\dots A_n)
 $$
where $A:=\pi(\alpha)$. We write $\Shs_M$ for the set of
explicitely Markovian sequences of $\Shs$.
\end{defi}

\begin{prop}
The restriction $\pi:\Shs_M\to\Sigma_M$ is a conjugacy.
\end{prop}

The proof of this proposition is the same as Proposition
\ref{prop-conjugacy-complete}

\begin{prop}
Let $\muh$ be a $\sigma$-invariant probability measure on $\Shs$
with $\muh(\Shs\setminus\Shs_M)=1$. Then
 $$
          h(\muh,\sigma) \leq \hcs(\Sigma).
 $$
\end{prop}

\begin{exam}
We prove that the above inequality can be reached and that the
entropy of the measure can exceed $\hc(\Sigma)<\hcs(\Sigma)$.

Let $\Sigma$ be the sofic subshift over the alphabet $\{0,1,2\}$
defined by the condition: between two zeroes there are an even
number of non-zero symbols. The Hofbauer's Markov diagram contains
the following vertices:
 \begin{enumerate}
    \item $\fol(0)$;
    \item $\fol(01)$,  $\fol(02)$;
    \item $\fol(011)$, $\fol(022)$;
    \item $\fol(1)$, $\fol(2)$.
 \end{enumerate}
The arrows are the following:
 \begin{itemize}
    \item $\fol(0)$ points to all the vertices of type 2;
    \item each vertex of type 2 points to all the vertices of type
3;
    \item each vertex of type 3 points to $\fol(0)$ as well as to
all the vertices of type 2;
    \item each vertex of type 4 points to all vertices of type 4
as well as to $\fol(0)$.
 \end{itemize}
It is easy to check:
\begin{itemize}
    \item $\hc(\Sigma)=0$ ($\Sigma$ is sofic);
    \item $\hcs(\Sigma)=\log 2$ ($\MCs(n)$ is the set of words $0\{1|2\}^{n-1}$;
    \item the non-explicitely Markovian sequences are exactly the
paths living on vertices of type 2 and 3 only.
\end{itemize}
\end{exam}

\begin{demof}{the proposition}
We fix $\mu$ as above and bound its entropy. Without losing
generality, we assume $\muh$ to be ergodic. We denote by $\hat{N}$
the set of sequences in $\Sh$ which are not explicitely Markovian.
We first claim that:
 $$
  \Sh\setminus \Sh_M = \hat{N}
 $$
up to a $\muh$-negligible set. Indeed, if $\alpha$ is explicitely
Markovian, then
 $$
     \alpha_0 = \fol (A_{-n}\dots A_0)
 $$
for all $n\geq0$ (where $A=\pi(A)$). Applying $\sigma_+(\cdot)\cap
[A_1]$ to each side of this inclusion
 $$
    \alpha_1 = \fol(A_{-n}\dots A_0A1).
 $$
Hence $\sigma^{-1}(\hat{N})\subset \hat{N}$. As in the proof of
the previous proposition, the claim follows.

\medbreak

Let $\mu=\pi_*\muh$. We remark that $h(\muh,\sigma)=h(\mu,\sigma)$
as $\pi:\Sh\to\Sigma$ is countable-to-one. We bound the entropy of
$\mu$ by bounding the minimal number of $n$-cylinders whose union
has measure $>1/2$ (see, e.g., \cite{DR}). Let $\eps>0$.

Let
 $$
   \ell(\alpha)=\min\{\ell\geq0:\exists B_{-\ell}\dots B_0
\text{ s.t. } \alpha_0=\fol(B_{-\ell}\dots B_0)\}.
 $$
$\ell(\alpha)<\infty$ everywhere. Hence, one can find $L_0<\infty$
such that $\{\alpha\in\Sh:\ell(\alpha)<L_0\}$ has positive
$\muh$-measure.

Let $K_0<\infty$ be such that $\#\MC(n)\leq
e^{(\hc(\Sigma)+\eps)n}$ for all $n\geq K_0$. We assume $K_0\geq
L_0/\eps$.

Let $n(\alpha)=\min\{k\geq K_0:\ell(\sigma^{-k}\alpha)< L_0\}$.
This is well-defined $\muh$-a.e. by ergodicity. There exists
$N_0<\infty$ such that $n(\cdot)>N_0$ on a set of $\muh$-measure
$<\eps/\log\#\A$.

By Birkhoff's ergodic theorem, there exist an integer $M_0<\infty$
and a measurable set $\hat G_0\subset\Sh$ with $\mu(G_0)>1/2$ such
that for all $\alpha\in \hat G_0$, all $n\geq M_0$,
 $$
   \frac1n\#\{ 0\leq k<n : n(\sigma^k\alpha)>N_0\} <
\frac{\eps}{\log\#\A}.
 $$
We may assume that $M_0\geq N_0 \log\#\A/\eps$.

It is easy to see that for any $n\geq M_0$, any $\alpha\in G_0$,
$A_0\dots A_{n-1}$ (recall $A=\pi(\alpha)$) can be decomposed
into:
 \begin{itemize}
  \item segments of the form $A_{m-k}\dots A_m$ with $k:=n(\sigma^m\alpha)<N_0$;
  \item an initial segment of length at most $N_0$;
  \item at most $\eps n/\log\#\A$ left-overs.
 \end{itemize}

Let us show that the segments of the first type are essentially
left constraints. Fix one such segment $w=A_{m-k}\dots A_m$ with
$k:=n(\sigma^m\alpha)$. $\alpha_{m-k}=\fol(B_{-p}\dots B_0)$ for
some finite word of length $p\leq L_0$ with $B_0=A_{m_k}$ of
course). $\sigma^{m-k}\alpha$ is not explicitely Markovian, i.e.,
 $$
    \alpha_m = \fol(B_{-p}\dots B_0A_{m-k+1}\dots A_m)
             \subsetneq \fol(A_{m-k}\dots A_m)
 $$
Hence, $\fol(B_{-p}\dots B_0)$ meets but does not include
$[A_{m-k}\dots A_m]$. Let $C\in[A_{m-k}\dots A_m]\setminus
\fol(B_{-p}\dots B_0)$ so that $B_{-p}\dots
B_{-1}C_0C_1\dots\notin\Sigma_+$. Hence, if we let $q$ be the
smallest integer such that $B_{-q}\dots
B_{-1}C_0C_1\dots\notin\Sigma_+$, we get that:
 $$
        B_{-q}\dots B_{-1}C_0\dots C_k\in\MCs(k+q)
 $$
Notice that $q\leq \eps n$.

 Thus, the number of $n$-cylinders meeting $\pi(\hat G_0)$
(which has $\mu$-measure $>1/2$) is bounded by:
 $$
    C^{2n/K_0}_n e^{(\hcs(\Sigma)+\eps)(1+\eps)n} \#\A^{N_0+\eps
n/\log\#\A}
       \leq e^{(\hcs(\Sigma)+4\eps)(1+\eps)n}.
 $$
As $\eps>0$ was arbitrary, this proves that $h(\mu,\sigma)\leq
\hcs(\Sigma)$.
\end{demof}

\end{document}